\documentclass[12pt]{article}
 \usepackage{latexsym}
 \usepackage{amssymb}
 \usepackage{color}
 \usepackage{graphicx}
\usepackage{epstopdf}
\usepackage{amsmath}    %  Added by Yannan
\usepackage{bm}                     % Define vectors, Added by Yannan
 \newtheorem{Theorem}{Theorem}

\newcommand{\A}{{\cal A}}

\newcommand{\uu}{{\bf u}}

\newcommand{\x}{{\bf x}}

\newcommand{\qed}{\nobreak \ifvmode \relax \else
      \ifdim\lastskip<1.5em \hskip-\lastskip
      \hskip1.5em plus0em minus0.5em \fi \nobreak
      \vrule height0.75em width0.5em depth0.25em\fi}

\begin{document}
\title{An Explicit SOS Decomposition of A Fourth Order Four Dimensional Hankel Tensor with A Symmetric Generating Vector}

\author{Yannan Chen\footnote{School of Mathematics and Statistics,
Zhengzhou University, Zhengzhou, China. E-mail: ynchen@zzu.edu.cn (Y. Chen)
    This author's work was supported by the
    National Natural Science Foundation of China (Grant No. 11401539) and
    the Development Foundation for Excellent Youth Scholars of
    Zhengzhou University (Grant No. 1421315070).}
    \quad Liqun Qi\footnote{Department of Applied
 Mathematics, The Hong Kong Polytechnic University, Hung Hom, Kowloon, Hong
 Kong.
 E-mail: maqilq@polyu.edu.hk (L. Qi).  This author's work was partially supported by the Hong Kong Research Grant Council (Grant No. PolyU
502111, 501212, 501913 and 15302114).}   \quad   Qun Wang \footnote{Department of Applied
 Mathematics, The Hong Kong Polytechnic University, Hung Hom, Kowloon, Hong
 Kong.   Email: wangqun876@gmail.com (Q. Wang).} }

\date{\today} \maketitle

%---------------------------------------------------------------------------------Abstract
\begin{abstract}
\noindent  %\vspace{3mm}
 In this note, we construct explicit SOS decomposition of A Fourth Order Four Dimensional Hankel Tensor with A Symmetric Generating Vector, at the critical value.   This is a supplementary note to Paper \cite{CQW}.

\noindent {\bf Key words:}\hspace{2mm} Hankel tensors, generating vectors, sum of squares, positive semi-definiteness.

\noindent {\bf AMS subject classifications (2010):}\hspace{2mm}
15A18; 15A69
  \vspace{3mm}

\end{abstract}

We construct an explicit SOS decomposition for the homogeneous polynomial $f(\x)$.
Let $$ \uu := \left(x_1^2, x_2^2, x_3^2, x_4^2, x_1x_3, x_2x_4,
      x_1x_2, x_3x_4, x_2x_3, x_1x_4\right)^\top, $$
be a basis of quadratics.
Then, $f(\x)$ is SOS if and only if there exists a PSD matrix $C$ such that $f(\x)=\uu^\top C\uu$.
Inspired by the Cholecky decomposition of $C$,
we define the parameterized SOS decomposition of $f(\x)$ as follows
$$ f(\x)=\sum_{k=1}^{10} q_k^2(\x), $$
where
\begin{eqnarray*}
  q_1(\x) &=& \alpha_{11}x_1^2+\alpha_{12}x_2^2+\alpha_{13}x_3^2+\alpha_{14}x_4^2+\alpha_{15}x_1x_3+\alpha_{16}x_2x_4+\alpha_{17}x_1x_2+\alpha_{18}x_3x_4 \\
          &&{}+\alpha_{19}x_2x_3+\alpha_{1,10}x_1x_4, \\
  q_2(\x) &=& \alpha_{22}x_2^2+\alpha_{23}x_3^2+\alpha_{24}x_4^2+\alpha_{25}x_1x_3+\alpha_{26}x_2x_4+\alpha_{27}x_1x_2+\alpha_{28}x_3x_4+\alpha_{29}x_2x_3 \\
          &&{}+\alpha_{2,10}x_1x_4, \\
  q_3(\x) &=& \alpha_{33}x_3^2+\alpha_{34}x_4^2+\alpha_{35}x_1x_3+\alpha_{36}x_2x_4
            +\alpha_{37}x_1x_2+\alpha_{38}x_3x_4+\alpha_{39}x_2x_3+\alpha_{3,10}x_1x_4, \\
  q_4(\x) &=& \alpha_{44}x_4^2+\alpha_{45}x_1x_3+\alpha_{46}x_2x_4
            +\alpha_{47}x_1x_2+\alpha_{48}x_3x_4+\alpha_{49}x_2x_3+\alpha_{4,10}x_1x_4, \\
  q_5(\x) &=& \alpha_{55}x_1x_3+\alpha_{56}x_2x_4
            +\alpha_{57}x_1x_2+\alpha_{58}x_3x_4+\alpha_{59}x_2x_3+\alpha_{5,10}x_1x_4, \\
  q_6(\x) &=& \alpha_{66}x_2x_4+\alpha_{67}x_1x_2+\alpha_{68}x_3x_4+\alpha_{69}x_2x_3+\alpha_{6,10}x_1x_4, \\
  q_7(\x) &=& \alpha_{77}x_1x_2+\alpha_{78}x_3x_4+\alpha_{79}x_2x_3+\alpha_{7,10}x_1x_4, \\
  q_8(\x) &=& \alpha_{88}x_3x_4+\alpha_{89}x_2x_3+\alpha_{8,10}x_1x_4, \\
  q_9(\x) &=& \alpha_{99}x_2x_3+\alpha_{9,10}x_1x_4, \\
  q_{10}(\x) &=& \alpha_{10,10}x_1x_4.
\end{eqnarray*}
The involved $55$ parameters $\alpha_{i,j}, (i \leq j),$ satisfies the following $35$ equality constraints
\begin{eqnarray}
     v_0 &=& \alpha_{11}^2, \label{ce2}\\
   4v_1 &=& 2\alpha_{11}\alpha_{17}, \\
   4v_2 &=& 2\alpha_{11}\alpha_{15}, \\
   6v_2 &=& 2\alpha_{11}\alpha_{12}+\sum_{k=1}^7\alpha_{k7}^2, \\
   4v_3 &=& 2\alpha_{11}\alpha_{1,10}, \\
  12v_3 &=& 2\alpha_{11}\alpha_{19}+2\sum_{k=1}^5\alpha_{k5}\alpha_{k7}, \\
   4v_3 &=& 2\sum_{k=1}^2\alpha_{k2}\alpha_{k7}, \\
      12 &=& 2\alpha_{11}\alpha_{16}+2\sum_{k=1}^7\alpha_{k7}\alpha_{k10}, \\
       6 &=& 2\alpha_{11}\alpha_{13}+\sum_{k=1}^5\alpha_{k5}^2, \\
      12 &=& 2\sum_{k=1}^2\alpha_{k2}\alpha_{k5}+2\sum_{k=1}^7\alpha_{k7}\alpha_{k9}, \\
       1 &=& \sum_{k=1}^2\alpha_{k2}^2, \\
  12v_5 &=& 2\alpha_{11}\alpha_{18}+2\sum_{k=1}^5\alpha_{k5}\alpha_{k,10}, \\
  12v_5 &=& 2\sum_{k=1}^2\alpha_{k2}\alpha_{k10}+2\sum_{k=1}^6\alpha_{k6}\alpha_{k7}, \\
  12v_5 &=& 2\sum_{k=1}^3\alpha_{k3}\alpha_{k7}+2\sum_{k=1}^5\alpha_{k5}\alpha_{k9}, \\
   4v_5 &=& 2\sum_{k=1}^2\alpha_{k2}\alpha_{k9}, \\
   6v_6 &=& 2\alpha_{11}\alpha_{14}+\sum_{k=1}^{10}\alpha_{k10}^2, \\
  24v_6 &=& 2\sum_{k=1}^5\alpha_{k5}\alpha_{k6}+2\sum_{k=1}^7\alpha_{k7}\alpha_{k8}+2\sum_{k=1}^9\alpha_{k9}\alpha_{k10}, \\
   4v_6 &=& 2\sum_{k=1}^3\alpha_{k3}\alpha_{k5},
\end{eqnarray}
\begin{eqnarray}
   4v_6 &=& 2\sum_{k=1}^2\alpha_{k2}\alpha_{k6}, \\
   6v_6 &=& 2\sum_{k=1}^2\alpha_{k2}\alpha_{k3}+\sum_{k=1}^9\alpha_{k9}^2, \\
  12v_5 &=& 2\sum_{k=1}^4\alpha_{k4}\alpha_{k7}+2\sum_{k=1}^6\alpha_{k6}\alpha_{k10}, \\
  12v_5 &=& 2\sum_{k=1}^3\alpha_{k3}\alpha_{k10}+2\sum_{k=1}^5\alpha_{k5}\alpha_{k8}, \\
  12v_5 &=& 2\sum_{k=1}^2\alpha_{k2}\alpha_{k8}+2\sum_{k=1}^6\alpha_{k6}\alpha_{k9}, \\
   4v_5 &=& 2\sum_{k=1}^3\alpha_{k3}\alpha_{k9}, \\
      12 &=& 2\sum_{k=1}^4\alpha_{k4}\alpha_{k5}+2\sum_{k=1}^8\alpha_{k8}\alpha_{k10}, \\
       6 &=& 2\sum_{k=1}^2\alpha_{k2}\alpha_{k4}+2\sum_{k=1}^6\alpha_{k6}^2, \\
      12 &=& 2\sum_{k=1}^3\alpha_{k3}\alpha_{k6}+2\sum_{k=1}^8\alpha_{k8}\alpha_{k9}, \\
       1 &=& \sum_{k=1}^3\alpha_{k3}^2, \\
   4v_3 &=& 2\sum_{k=1}^4\alpha_{k4}\alpha_{k10}, \\
  12v_3 &=& 2\sum_{k=1}^4\alpha_{k4}\alpha_{k9}+2\sum_{k=1}^6\alpha_{k6}\alpha_{k8}, \\
   4v_3 &=& 2\sum_{k=1}^3\alpha_{k3}\alpha_{k8}, \\
   4v_2 &=& 2\sum_{k=1}^4\alpha_{k4}\alpha_{k6}, \\
   6v_2 &=& 2\sum_{k=1}^3\alpha_{k3}\alpha_{k4}+2\sum_{k=1}^8\alpha_{k8}^2, \\
   4v_1 &=& 2\sum_{k=1}^4\alpha_{k4}\alpha_{k8}, \\
    v_0 &=& \sum_{k=1}^4\alpha_{k4}^2.  \label{ce3}
\end{eqnarray}

Giving fixed $v_2,v_6,v_1,v_3,v_5$, we denote $M_1(v_2,v_6,v_1,v_3,v_5)$ as the minimum of
the following optimization problems:
\begin{equation}\label{ce1}
\begin{aligned}
  M_1(v_2,v_6,v_1,v_3,v_5) := \min ~&~ v_0 \\
    \mathrm{s.t.} ~&~ \text{equality constraints (\ref{ce2})-(\ref{ce3}).}
\end{aligned}
\end{equation}

By simple algebraic derivation, we have the following theorem.

\begin{Theorem} \label{t12}
Suppose that the assumptions (2) and (3) in \cite{CQW} hold.  Then if $v_0 \ge M_1(v_2,v_6,v_1,v_3,v_5)$,
$\A$ is SOS and
$$f(\x) = \sum_{k=1}^{10} q_k(\x)^2.$$
Thus, $M_1(v_2,v_6,v_1,v_3,v_5) \geq M_0(v_2,v_6,v_1,v_3,v_5)$.
\end{Theorem}

%%%%%%%%%%%%%%%%%%%%%%%%%%%%%%%%%%%%%%%%%%%%%%%%%%%%%%%%%%%%%%%%%%%%%%%%%%%%%%%%%%%%%%%%%%%%%%%%%%%%%%%%%%%%%%%%%%%%%%%%%%%%%%%%%%%%%%%%%%%%%%%%%%%%%%%%%%%%

\end{document}